\documentclass[graybox]{svmult}

\usepackage{mathptmx}     
\usepackage{helvet}       
\usepackage{courier}      
\usepackage{type1cm}      
                          
\usepackage{makeidx}   
\usepackage{graphicx}  
                       
\usepackage{multicol}
\usepackage[bottom]{footmisc}
\usepackage{amsmath,amssymb}

\newcommand{\R}{\mathbb{R}}

\makeindex      

\begin{document}

\title*{A New Mathematical Model for the Efficiency Calculation\thanks{This is a preprint 
of a paper accepted for publication 11-March-2019 as a book chapter 
in ``Advances in Intelligent Systems and Computing'' 
(\url{https://www.springer.com/series/11156}), Springer.}}

\author{An\'ibal Galindro \and Micael Santos \and Delfim F. M. Torres \and Ana Marta-Costa}

\authorrunning{A. Galindro \and M. Santos \and D.F.M. Torres \and A. Marta-Costa} 

\institute{An\'ibal Galindro 
\at Centre for Transdisciplinary Development Studies,
University of Tr\'as-os-Montes and Alto Douro,
Polo II--ECHS, Quinta de Prados, 5000-801 Vila Real, Portugal.\\ 
\email{anibalg@utad.pt}
\and 
Micael Santos 
\at Centre for Transdisciplinary Development Studies,
University of Tr\'as-os-Montes and Alto Douro,
Polo II--ECHS, Quinta de Prados, 5000-801 Vila Real, Portugal.\\ 
\email{micaels@utad.pt}
\and 
Delfim F. M. Torres  (corresponding author)
\at Center for Research and Development in Mathematics and Applications (CIDMA),\\ 
Department of Mathematics, University of Aveiro, 3810-193 Aveiro, Portugal.\\ 
\email{delfim@ua.pt}
\and Ana Marta-Costa 
\at 
Centre for Transdisciplinary Development Studies,
University of Tr\'as-os-Montes and Alto Douro,
Polo II--ECHS, Quinta de Prados, 5000-801 Vila Real, Portugal.\\ 
\email{amarta@utad.pt}}

\maketitle


\abstract{During the past sixty years, a lot of effort has been made regarding 
the productive efficiency. Such endeavours provided an extensive bibliography 
on this subject, culminating in two main methods, named the Stochastic Frontier 
Analysis (parametric) and Data Envelopment Analysis (non-parametric). 
The literature states this methodology also as the benchmark approach, 
since the techniques compare the sample upon a chosen ``more-efficient'' reference. 
This article intends to disrupt such premise, suggesting a mathematical model 
that relies on the optimal input combination, provided by a differential 
equation system instead of an observable sample. A numerical example 
is given, illustrating the application of our model's features.}

\medskip

\noindent \textbf{Keywords:} 
frontier estimation models,
efficiency analysis,
technical inefficiency,
from data to differential equations.

\smallskip

\noindent \textbf{2010 Mathematics Subject Classification:} 90B50; 90C08. 


\section{Introduction}
\label{sec:1}

The word \emph{efficiency} may acquire several meanings, but in economics 
such word is closely related to the general premise of the field, since the 
available resources are limited and how should we use them to attain 
the maximum level of output or the maximum utility is a natural question.
The productivity efficiency literature began in 1957 with Farrell's work 
and is a rising theme applied to several sectors. Farrell \cite{farrell} 
defined two concepts of efficiency: Technical Efficiency (TE) and Allocative Efficiency (AE). 
The first is evident when a certain level of inputs is given, the Decision Making Unit (DMU) 
is able to produce the maximum level of output or, fixing a certain level of output, 
the DMU is able to minimize the level of input \cite{bravo,fleming}. The AE reflects 
the ability of a firm to use the inputs in their optimal proportions (given their respective prices), 
to minimize the cost or maximize the revenue \cite{bravo,aparicio}. Another type of efficiency 
is the scale efficiency, which tells us if a DMU is operating on an optimal scale \cite{anang}.
Over time, several methodologies have been developed to estimate the productive efficiency. 
Those methodologies can be categorized mainly as parametric or non-parametric. The most 
or the main methods used in each category are the Stochastic Frontier Analysis (SFA) 
and Data Envelopment Analysis, respectively \cite{bravo,djokoto1,djokoto2,mareth,thiam}. 
However, with the development of the methodologies, such separation is not so clear nowadays. 
For example, Johnson and Kuosmanen \cite{johnson} have developed a semi-nonparametric 
one-stage DEA and many efforts have been made to develop non-parametric or semiparametric 
SFA, as documented on the work of Kumbhakar, Parmeter, and Zelenyuk \cite{kumbhakar}. 
However, all these methodologies are benchmarking approaches, in other words, they compare 
the productive efficiency of a DMU against a reference performance (the most-efficient DMU's). 
In contrast, the goal of this work is to create a new methodology that does not follow 
this benchmark premise and looks beyond the self-proclaimed better or more efficient observation. 
Our differential equation model relies on overall optimality per delivered output, 
which endures a different interpretation of the problem.


\section{Literature Review}
\label{sec:2}

The productive efficiency methodologies lay on the benchmarking approach, since all 
the techniques rely on a continuous and systematic process of comparing a certain 
chosen sample upon a reference (benchmark) performance \cite{jamasb,khetrapal}. 
The reference performance is normally the ``best-practice'', i.e., the methodologies 
identify the most efficient DMUs to build a frontier and then compare the least 
efficient DMUs against this frontier. The greater the distance from a DMU performance 
to the best-practice, the greater its level of inefficiency. The benchmark approach 
in the productive efficiency can be categorized mainly in two groups: non-parametric 
and parametric. Non-parametric techniques are normally based on programming techniques 
and do not require a production, cost, or profit function, calculating the relation 
of the inputs with the outputs without an econometric estimation \cite{bravo,thiam,khetrapal}. 
Within the non-parametric approaches, we have the DEA and the Free Disposal Hull (FDH).
Usually, the productive efficiency methods assume the convexity of the production set. 
However, FDH was created by Deprins, Simar, and Tulkens \cite{deprins}, who proposed 
the elimination of the convexity assumption \cite{leleu}. FDH is a variation of DEA, 
based on a programing technique, and the frontier estimated by this method may coincide 
or be below the DEA frontier, so the efficiency scores estimated by FDH tend to be higher 
than those estimated by DEA \cite{tulkens}. The DEA is the most used non-parametric TE model. 
It is based on mathematical programming techniques and does not require the specification 
of a functional form for the technology. It should be noted that the majority of 
the non-parametric methods are deterministic, therefore they do not allow any random noises 
or measurement errors \cite{bravo,thiam}. Another characteristic of the DEA methods 
is their potential sensitivity of efficiency scores to the number of observations, 
as well as to the dimensionality of the frontier and to the number of outputs and inputs 
\cite{bravo,thiam,ramanathan}. However, some developments have been made in the 
traditional DEA. For example, Simar and Wilson developed a stochastic DEA using  
bootstrapping techniques in \cite{simar1,simar2}. Johnson and Kuosmanen \cite{johnson} 
developed a semi-nonparametric one-stage DEA based on the critics of the two-stage DEA 
and on the work of Banker and Natarajan \cite{banker}, who incorporate a noise term 
that has a truncated distribution. SFA normally is a parametric and stochastic method, 
which was introduced by Aigner, Lovell, and Schmidt \cite{aigner} and Meeusen 
and van Den Broeck \cite{meeusen}. This method imposes more structure on the shape 
of the frontier by specifying a functional form for the production function, 
such as the Cobb--Douglas or the Translog form. Moreover, the SFA was developed 
to allow random errors, but neither the random error nor inefficiency can be observed, 
so separating them requires an assumption on the distribution of the efficiencies 
scores and on the random error. Although the SFA is traditionally parametric, 
some developments have given it some degree of convergence with non-parametric models, 
which are referred as non-parametric or semiparametric SFA 
\cite{kumbhakar,banker2,fan,kuosmanen,noh,parmeter}. Essentially, Fan et al. \cite{fan} 
and Kneip and Simar \cite{kneip} provide the baseline studies for these methodologies. 
More recently, the advances in stochastic frontier models were documented 
in Kumbhakar et al. \cite{kumbhakar}.


\section{A New Theoretical Model}
\label{sec:3}

The main goal of this article is to overcome the limitations of the benchmark approach, 
extensively used in the previously presented bibliography, appended to several efficiency 
calculation methods. The generalized method of obtaining a sample (at least one)
based totally on efficient benchmark does not guarantee a proper or correct solution. 
Therefore, the subsequent process of evaluating the efficiency/inefficiency level of 
the other observations might be biased. We propose a differential equation based method 
to solve a single output and multiple input efficiency problem. Let $m$ be the sample 
size containing the data points for $n$ input variables $X=\{x_1,x_2,x_3,\ldots,x_n\}$ 
and $Y$ the output associated to each data point such that $Y=\{y_1,y_2,y_3,\ldots,y_m\}$. 
The explicit matricial form is stated on Equation~\eqref{matrix}:
\begin{equation}
\label{matrix}
\begin{bmatrix} 
{y^1} \\ {y^2} \\ \vdots \\ {y^m} 
\end{bmatrix}
=
\begin{bmatrix} 
{x_1}^{1}& {x_2}^{1} & \dots & {x_n}^{1} \\
{x_1}^{2} & {x_2}^{2} & \dots & {x_n}^{2} \\
\vdots & \vdots & \ddots & \vdots \\
{x_1}^{m}      & {x_2}^{m} & \dots & {x_n}^{m}
\end{bmatrix}.
\end{equation}
The main idea is to interpret the input/output dynamics as a differential equation 
system that retrieves the optimal input combination $X^*$ per each output level 
$Y\in{\R^+}$. Therefore, the system compiles $n$ differential equations, one per each 
input variable, and it is defined by Equation~\eqref{diff}, where 
$\bold{A_i}=\{{\beta_i}^{0},{\beta_i}^{1},\ldots,{\beta_i}^{n}\}$, 
$i=1,2,\ldots,n$, represents the whole inner input trade-off combinations that forces 
optimal input allocation per each output level. Such parameters can be numerically 
estimated for both linear and non-linear differential equations using methods 
developed by Ramsay~\cite{ramsay,ramsay2}:   
\begin{eqnarray}
\label{diff}
\begin{cases}
\displaystyle \frac{dx_1(y)}{dy}
=\beta_1^0 + \beta_1^1 x_1(y) + \beta_1^2 x_2(y) + \cdots + \beta_1^n x_n(y) ,\\[0.3cm]
\displaystyle \frac{dx_2(y)}{dy}
=\beta_2^0 + \beta_2^1 x_1(y) + \beta_2^2 x_2(y) + \cdots + \beta_2^n x_n(y) ,\\[0.3cm]
\quad \vdots \\[0.3cm]
\displaystyle \frac{dx_n(y)}{dy}
=\beta_n^0 + \beta_n^1 x_1(y) + \beta_n^2 x_2(y) + \cdots + \beta_n^n x_n(y).
\end{cases}
\end{eqnarray}
We settle that the set $X$ of input variables is expressed in similar values 
(Euro, for example), to avoid the concern of different weights or costs in the 
model development. The optimality of $A^*$ does not guarantee the feasibility 
of the process with non-negative inputs across every output value. 
Therefore, the condition on Equation \eqref{rest} should be met:
\begin{equation}
\label{rest}
0 \leq x_1,x_2,\ldots,x_n \  \forall  \ y.
\end{equation}
The model initial conditions should follow Equation \eqref{IC}, 
where the null output is generated by zero inputs:
\begin{equation}
\label{IC}
x_1(0)=0, \quad x_2(0)=0, \quad \ldots, \quad x_n(0)=0.
\end{equation}
Nonetheless, we intend to also obtain the optimal output level $y^*$ 
for the given problem. Assuming increasing returns to scale until 
some single point (which would be $y^*$) and decreasing returns to scale 
afterwards, per each $y^d > y^*$, the point $y^*$ guarantees profit 
maximization and solution singleness. The profit function given by Equation~\eqref{Opt} 
acquaints that the resulting output is sold at an arbitrary price $c$. Each level 
of output $y$ is bundled with an aggregate variation rate of every input 
on Equation \eqref{Opt}. The integration extracts the numerical input 
for each level of production $y$. The maximizing profit stretched 
by the optimal $y^*$, is obtained when the derivative of $J(y)$, 
\begin{equation}
\label{Opt}
J(y) = cy - \int_0^{y} \frac{dx_i(y)}{dy},
\end{equation}
equals zero:
\begin{equation}
\label{Opt1}
\frac{dJ(y)}{dy} = 0.
\end{equation}
The newly discovered $y^*$ eases straightforwardly the computation 
of the ideal input set $X^*=\{x_1^*,x_2^*,x_3^*,\ldots,x_n^*\}$, 
where each real value can be obtained by the integration of Equation \eqref{int}:
\begin{equation}
\label{int}
x_i^*=\int_0^{y^*} \frac{dx_i(y)}{dy}.
\end{equation}
The optimal $X^*$ can be interpreted as a point in a $n$-dimensional space, 
each non-optimal observation $X_j^{no}=\{x_{1j}^{no},x_{2j}^{no},x_{3j}^{no},
\ldots,x_{nj}^{no}\} \in \R^n$, $j=1,2,3,\ldots,m$, from the sample $m$
outlies a distance vector from the optimal solution on Equation \eqref{vect}:
\begin{equation}
\label{vect}
d_j = \sqrt{(x_{1j}^{no}-x_1^*)^2+(x_{2j}^{no}-x_2^*)^2
+ \cdots + (x_{nj}^{no}-x_n^*)^2}.    
\end{equation}
Let $\Omega$ be the set of feasible input combinations such that 
$f:\Omega \subset \R^{n} \rightarrow \R$ is a function with domain $\Omega$ 
with values in $\R$ (the output level). Analogously, we settle a subset function 
using a restricted domain $\Omega' \subseteq \Omega$, which contains the feasible 
data points that obey a certain restriction. A given point from $\Omega'$ 
is extracted from the $\Omega$ set only if $\forall \Omega'$, $0 \leq \frac{dy}{dx_i}$, 
and the condition on Equation \eqref{par} is met:   
\begin{equation}
\label{par}
\sum_{i=0}^{n} \frac{dy}{dx_i} > \epsilon.  
\end{equation}  
Naturally, the selection of $\epsilon$ acquaints a certain level of parsimony, 
since larger values may absorb and exclude valid and empirical observations. 
On the other hand, choosing an $\epsilon$ that is too small may expand unnecessarily 
the feasible region $\Omega'$. Considering each element from the subset $\Omega'$, 
obtained with a certain level of $\epsilon$, it is possible to generate 
an analogous set $\Omega''$ that contains the distance given by 
Equation \eqref{vect} of each feasible input combination from $\Omega'$. 
Let $X_w$ be the input combination from $\Omega'$ that obtains the maximum value 
from the $\Omega''$ set ($d_w$), which is collinear with the observation that 
we intend to attain information about their efficiency. Let it be $X_m^{no}$. 
Literally, $X_m^{w}$ work as the most inefficient points and as a general 
reference for the efficiency calculations. Being 
$X_j^{no}=\{x_{1j}^{no},x_{2j}^{no},x_{3j}^{no},
\ldots,x_{nj}^{no}\} \in \R^n$, $j=1,2,\ldots,m$, the sampled non-optimal data points, 
with distances to $X^*$ given by Equation \eqref{vect}, 
their efficiency level is given by Equation \eqref{eff}:
\begin{equation}
\label{eff}
I_d = 1 - \frac{d_j}{d_w}, 
\end{equation}
for all $d_j \leq d_w$. Therefore, the efficiency level of each sampled 
observation $I_d$ in bounded between 0 and 1. Since the weights or costs are 
acquainted, since the beginning, within the input variables, the distances 
on the generated $n$-dimensional space provide the same effort throughout the efficient point. 
Let $X_a$ and $X_b$ be to sampled data points such that $I_a=I_b$ is the equidistant 
measurement provided in the $n$-dimensional space embedded with the pre-accounted 
weights, ensuring that the effort through the optimal point (efficiency improvement) 
of both samples $X_a$ and $X_b$ is the same, validating the fairness when we measure their efficiency levels. 
Nonetheless, we can stretch forward our model to acquaint singularly technical and allocative efficiency, 
instead of the previous general effort throughout the optimal point $y^*$. If we reduce the domain of $f$ 
even further, and assume subsets $\Omega^{m}$, $m \in \mathbb{N}$, we get $m$ functions 
$f_m$ such that $f_m:\Omega^{m} \subset \R^{n} \rightarrow \mathbb{R}$ satisfies
$f_m(x) = y_{m}$ for all $x \in \Omega^{m}$. Each given function $f_m$ 
with domain $\Omega^{m}$ acquaints a fixed output level $y^m$. Therefore, the resulting optimal points, 
considering that there is only an optimal set of inputs, is now given by $X^{m*}$ bundled with $y_{m*}$, 
extracted from the optimal movement of the differential equations on Equation~\eqref{diff}. 
Now, each sampled observation lays on a different $f_m$ for each output level $y_m$.
Assuming a sampled non-optimal input combination $X_m^{no}$, entangled with a certain 
level of output $y_m$, the productive efficiency is now given by the ratio of the distance 
to the optimal point $X^{m*}$ (Equation~\eqref{eff}) alongside the distance from the worst 
feasible point $X_m^{w}$ considered to $X^{m*}$. It is also possible to analogously shift 
the formulation to attain information about the technical efficiency of each observation. 
Instead of working with a subset that compiles similar levels of output, such subset 
should only feature similar levels of input given by the explicit sum of the input combination.  

Summarizing, we have just obtained the following result.

\begin{theorem}
Let $\left\{{\beta_i}^{0},{\beta_i}^{1},\ldots,{\beta_i}^{n}\right\} \in \R$, 
$i=1,2,\ldots,n$, compile the minimal input differential equation 
system \eqref{diff} per each $y$. Settle the initial conditions as \eqref{IC}, 
while assuring that the non-negativity condition \eqref{rest} holds.
A feasible $\Omega$ region is obtained with \eqref{par}.
Using the optimal $X^*$ per level of output $y$, the frontier points 
of $\Omega$ yield the worst efficiency levels $X_m^{w}$.
Sampling an observation $X_j$ with the associated collinear most 
inefficient point $X_w$, the distance vector is calculated according 
to \eqref{vect} using the optimal point $X^*$.
Finally, the efficiency levels are given by \eqref{eff}.
\end{theorem}


\section{An Illustrative Numerical Application}
\label{subsec:1}

\begin{figure}[t]
\sidecaption
\includegraphics[scale=.55]{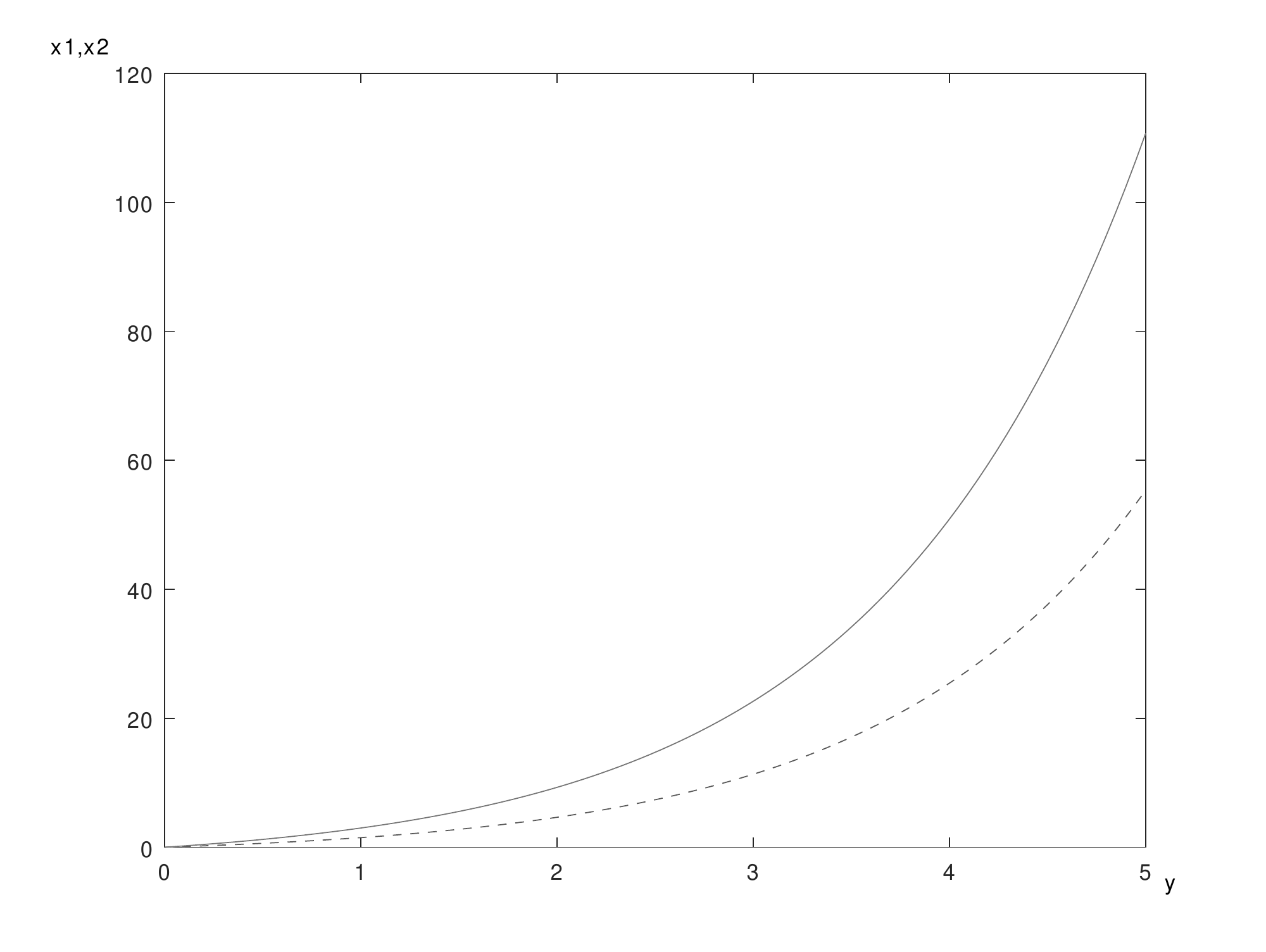}
\caption{The dashed/continuous line represents the optimal input level  
$x_1$/$x_2$ per each level of output $y$, respectively.}
\label{fig:1}     
\end{figure}

\begin{figure}[t]
\sidecaption
\includegraphics[scale=.65]{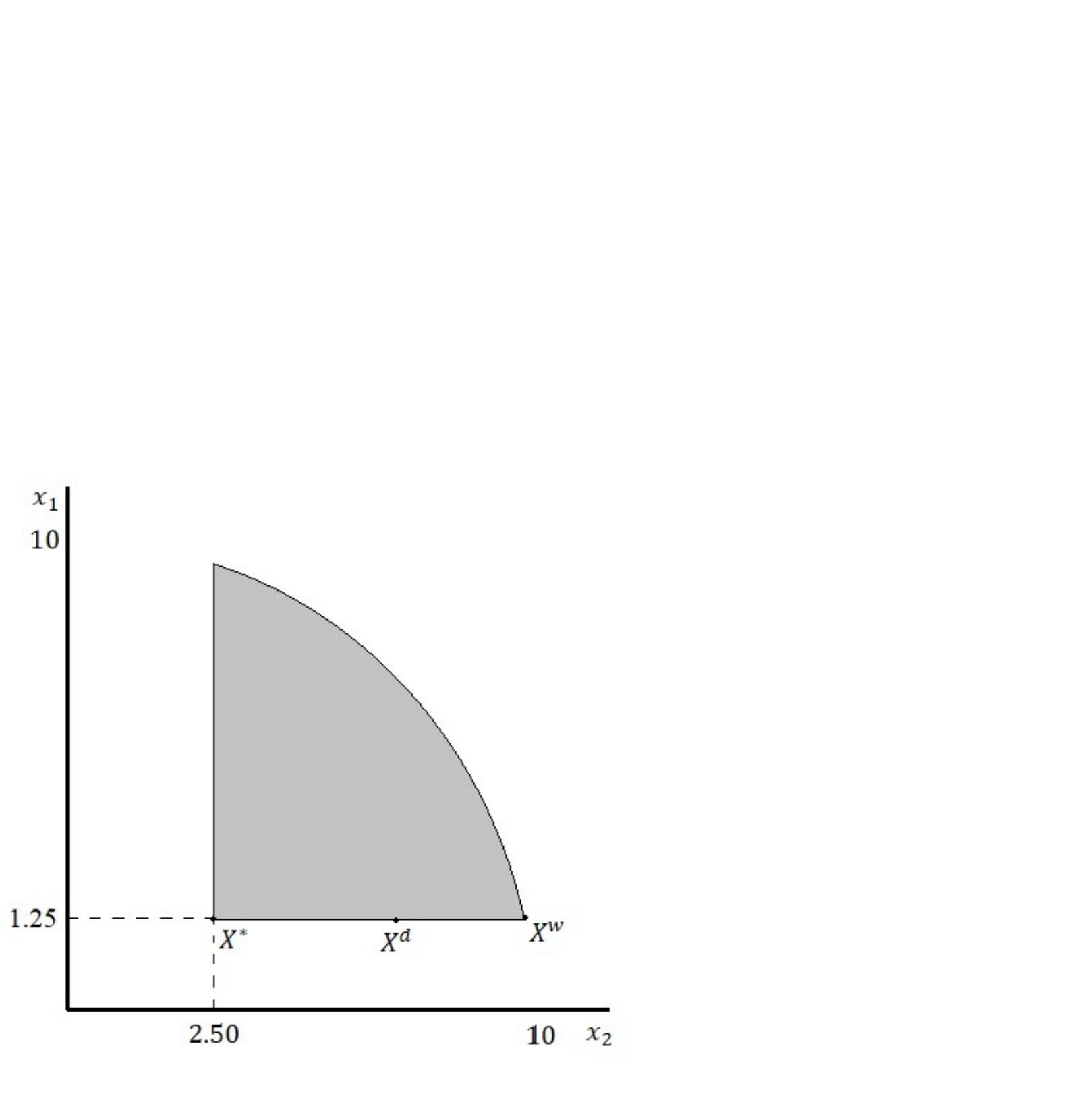}
\caption{The grey area represents the admissible input combinations $\Omega$ 
while $X^*$, $X^d$, and $X^w$ represent the optimal, sampled and worst collinear 
input levels, respectively. The output is settled in $y=1$.}
\label{fig:2}     
\end{figure}

In order to illustrate our method, we give here a simple numerical example. 
Following the differential equation system premise on Equation \eqref{diff}, 
a two input/single output system is created according to Equation \eqref{diff2}:
\begin{eqnarray}
\label{diff2}
\begin{cases}
\displaystyle \frac{dx_1(y)}{dy}
=1 + 0.25 x_1(y) + 0.25 x_2(y),\\[0.3cm]
\displaystyle \frac{dx_2(y)}{dy}
=2 + 0.50 x_1(y) + 0.50 x_2(y).
\end{cases}
\end{eqnarray}        
The aforementioned system describes a limited capacity productive system 
with decreasing returns to scale. Indeed, the input needed to increment the output 
values grows exponentially to bulky and unrealistic combinations. Such system does
not violate the non-negativity condition on Equation \eqref{rest} for any $y>0$, 
considering also the initial conditions stated on Equation \eqref{IC} 
for both $x_1$ and $x_2$. The Figure~\ref{fig:1} displays the differential 
equation development considering different levels of output $y$. The $c$ value 
is chosen such that the Equations \eqref{Opt} and \eqref{Opt1} of the maximizing 
profit hold. To obtain a specific example, we chose the optimal input combination 
for $y^1=1$, which using the integral on Equation \eqref{int} 
yields $x_1=1.25$ and $x_2=2.5$.  
The next step is to settle the acceptable domain $\Omega$, 
which does not violate the Equation \eqref{par}, for an almost 
neglectable $\epsilon$. Such condition of $\Omega$ 
lay on the Equation \eqref{Space} restrictions:
\begin{equation}
\label{Space}
\begin{split}
\displaystyle x_1 &\geq 1.25 ,\\[0.3cm]
\displaystyle x_2 &\geq 2.50 ,\\[0.3cm]
\displaystyle x_1^2 + x_2^2 &\leq 100.
\end{split}
\end{equation}     
With the admissible space settled, it is now possible to proceed with concrete 
computations. Let us take the example of a company producing $y=1$ with input 
combinations settled on the point $X^d=[x_1,x_2]=[1.25, 6.21]$ 
with the worst directional choice $X_w$ defined 
as $X^w=[x_1^w,x_2^w]=[1.25, 9.92]$ (see Figure~\ref{fig:2}). 
The distance vectors $d_j$ and $d_w$ are easily obtained from Equation \eqref{vect}, 
yielding the values $3.71$ and $7.42$, respectively, with $d_w$ retrieving 
the worst case scenario considering that the chosen point belongs to it's vector. 
The efficiency level of our example $I_d$ is delivered on Equation \eqref{result}: 
\begin{equation}
\label{result}
I_d = 1 - \frac{3.71}{7.42} = 0.5. 
\end{equation}
Since the point is equidistant to the optimal point, 
as the worst case scenario we can state that both the inneficiency 
and efficiency levels are $0.5$. 


\section{Conclusions}
\label{sec:4}

We presented a literature review of the state-of-the-art efficiency 
methods on benchmark settling techniques based on existing data. 
In order to overcome the sensibility of the existing data, 
we proposed a new non-sample method based on differential equations 
that mimics the target single-output/multiple-output productive system. 
A few assumptions were theorized in order to obtain a differential equation 
system that retrieves the optimal input level per each output. 
Such generalizations allow us to settle an optimality condition to obtain 
the utopic input/output level that maximizes the overall profit. Using this 
result as a reference, the sample efficiency level is calculated in a generated 
sub-space that acquaints both feasible and rational data points in $\R^n$. 
Since the inputs already acquaint their inner weights, the $n$-dimensional space 
(considering $n$ inputs) distance measuring, assures proportionality among the sample. 
The same premise can be followed to attain productive and technical efficiency/inefficiency 
levels, considering alternative endeavours of the selected model subset. Our theoretical 
approach relies on unbiased calculations of the coefficients $\beta_i^j$ 
on Equation \eqref{diff}. Introducing error with real data may demote 
the global optimum output $y^*$ to an estimated one. Nonetheless, 
as long as the weights or costs are directly observable, 
the sample efficiency level should withstand even though 
the optimal point is biased. The presented numerical simulation displays 
a simple two input and one output example, where the efficiency level of a sampled 
company was calculated. Nonetheless, data based numerical simulations 
can also be computationally expensive when the number of inputs grows. 
Our model can be expanded to multiple output approaches using the same premisses.


\begin{acknowledgement}
This work was supported by the R\&D Project INNOVINE \& WINE Vineyard and Wine Innovation
Platform Operation NORTE-01-0145-FEDER-000038, co-funded by the European and Structural
Investment Funds (FEDER) and by Norte 2020 (Programa Operacional Regional do Norte 2014/2020).
Torres was supported by FCT through CIDMA, project UID/MAT/04106/2019.
\end{acknowledgement}



\end{document}